%% file: paperV3.tex
\documentclass{article}
\input{command}
\renewcommand{\remi}[1]{{ #1}}
\newcommand{\Remi}[1]{{  #1}}
\begin{document}
\title{Evaluating a distance function}
\author{R. Abgrall}
\date{\today}
\maketitle
\begin{abstract}
Computing the distance function to some surface or line is a problem that occurs very frequently. There are several ways of computing a relevant approximation of this function, using for example technique originating from the approximation of Hamilton Jacobi problems, or the fast sweeping method. Here  we make a link with some elliptic problem and propose a very fast way to approximate the distance function.
\end{abstract}
\section{Introduction}
In many cases, one has to evaluate the distance function  to a surface $\Gamma_D$ \remi{which is part of the boundary of an open set $\Omega\in \R^d$}. An example in fluid mechanics is that of turbulence modelling: in some models, one of the parameters in the evaluation of the turbulent viscosity is  the distance to the airfoil.  Other examples can be found in medical image processing, surface reconstruction, etc. To evaluate the distance function, there are many ways. One technique is to numerically evaluate the viscosity solution of the Eikonal equation, 
\begin{equation}
\label{eq:1}
\begin{split}
\Vert \nabla u\Vert -1=0 & \text{~}\bx \in\remi{ \Omega}\\
u(\bx)=0 & \text{ if and only if } \bx\in \Gamma_D\\
u(\bx)\geq 0.
\end{split}
\end{equation}
In this formulation, only one boundary condition is prescribed, the Dirichlet one on  $\Gamma_D$, nothing is said for the other parts of $\partial \Omega$.  Numerical techniques for this can be found in \cite{shu,abgrall}. This can be improved by using a fast marching technique in the spirit of \cite{sethian}, or technique \remi{coming} from computer sciences.  From time to time, one can see in the literature methods that compute an approximation of the distance function as the solution of some Laplace equation, an example can be found in \cite{nishi} and the reference therein. \remi{Other references, where elliptic problems are considered, can be found in \cite{kuzmin,adams}. They solve \eqref{eq:1} by minimizing $\big (\Vert \nabla u\Vert-1)^2$ using a variational formulation with continuous or discontinuous finite element that needs to be penalized  to enforce the Dirichlet boundary conditions. It is interesting to notice that their formulation is relatively close to ours, using completely different paths.}

 There is no obvious links between \eqref{eq:1}, which is of hyperbolic nature, and an elliptic problem. The purpose of this small note is to provide a link, via the Hopf-Cole transform, and to provide a very fast algorithm { (compared to explicit algorithms for computing the viscosity solution of \eqref{eq:1} to \Remi{evaluate} the solution of \eqref{eq:1})}, at least if one accept a small viscosity term (which will nevertheless be present in any numerical method of PDE origin), and discretisation errors.

The format of this note is the following. First, we recall the Hopf-Cole transform and show how it can be applied to the steady Eikonal equation. This leads to an elliptic problem on a function constructed from the distance. We discuss the boundary conditions for this problem, and in particular for the part of $\partial \Omega$ which is not $\Gamma_D$. Then we provide a numerical method, and show its behaviour on unstructured triangular meshes.

\section{The problem}
We want to solve  the following problem:
Let $\Omega\subset \R^n$ be open, and we set $\partial \Omega=\Gamma_D\cup \Gamma_S$, $\Gamma_D\cap \Gamma_S$ of empty interior.  We want to compute the distance function to $\Gamma_D$. We \remi{ consider} the problem of finding the viscosity solution of 
\begin{equation}
\label{HJ}
\begin{split}
\Vert \nabla u \Vert -1=0 & \quad \bx\in \Omega\\
u=0 & \quad\bx \in \Gamma_D\\
u=+\infty & \quad\bx \in \Gamma_S
\end{split}
\end{equation}
Of course $\Gamma_S$ can be empty, but $\Gamma_D$ is never empty by assumption. On $\Gamma_S$ we have set Soner type boundary condition, see \cite{barles} for example. 
For the sake of completeness, we recall the notion of viscosity solution for \eqref{HJ}:
Let $\varphi\in C^1(\overline{\Omega})$, and $\bx_0$ a point where $u-\varphi$ reaches a local minimum: \eqref{HJ} means that if $\bx_0\in \Omega$, 
$\Vert \nabla \varphi(\bx_0) \Vert -1 \geq 0,$ and if $\bx_0\in \Gamma_D$
$\min\big (\Vert \nabla \varphi(\bx_0) \Vert -1, u(x_0)\big )\geq 0$  while, if $\bx_0\in \Gamma_S$
$\Vert \nabla \varphi(\bx_0) \Vert -1\geq 0.$
If $\bx_0$ is a minimum of $u-\varphi$, we get:
 if $\bx_0\in \Omega$, 
 $\Vert \nabla \varphi(\bx_0) \Vert -1 \leq 0,$ if $\bx_0\in \Gamma_D$
$\max\big (\Vert \nabla \varphi(\bx_0) \Vert -1, u(x_0)\big )\leq 0$  and
there is no condition on $\Gamma_S$.

Here we propose a method where we solve a viscous regularisation of \eqref{HJ}, or more precisely  of
\begin{equation}
\label{HJ2}
\begin{split}
\Vert \nabla u \Vert^2 -1=0 & \quad \bx\in \Omega\\
u=0 & \quad\bx \in \Gamma_D\\
u=+\infty & \quad\bx \in \Gamma_S
\end{split}
\end{equation}
since the two problems have the same solutions, 
 we consider, for $\nu>0$, the problem (in the viscosity sense, see \cite{barles} for second order problems)
\begin{equation}
\label{viscous:HJ2}
\begin{split}
\Vert \nabla u \Vert^2 -1=\nu \Delta u & \quad \bx\in \Omega\\
u=0 & \quad\bx \in \Gamma_D\\
u=+\infty & \quad\bx \in \Gamma_S
\end{split}
\end{equation}

\section{Rewriting the problem}

If, instead of looking at the steady problem \eqref{viscous:HJ2}, we consider the unsteady one,
$$\dpar{u}{t}+\Vert \nabla u\Vert^2 -1=\nu \Delta u$$
with the same initial and boundary conditions, this \remi{is} "almost" the \remi{viscous}  Burgers equation, 
$$\dpar{u}{t}+\frac{1}{2}\Vert \nabla u\Vert^2 =\nu \Delta u.$$
for which  it is very well known  that we can transform it  into the heat equation by using the Hopf-Cole transform,
\begin{equation}\label{eq:3:3}u(\bx, t)=-2\nu\log\big (\varphi(\bx, t)\big ) . \end{equation}
The proof is classical (though done in one dimension in most textbooks), but we nevertheless repeat it.

Our notations will be: $\nabla u$ represents the first derivative of the function $u$: for any $h\in \R^d$
$$ u(\bx+h)=u(\bx)+ \nabla u(\bx)\cdot h+o(h), $$
and $D^2u$ represents the second derivative (the Hessian) of $u$:
$$\remi{\nabla u(\bx+h)=\nabla u (\bx)+D^2u(\bx)\cdot h +o(h)}.$$

With this in mind, we have, from \eqref{eq:3:3}, that
$$\nabla u(\bx,t)=-2\nu\dfrac{\nabla\varphi(\bx,t)}{\varphi(\bx,t)}, \qquad \dpar{u}{t}=-2\nu\dfrac{\dpar{\varphi(\bx,t)}{t}}{\varphi(\bx,t)}$$
and 
$$D^2u(\bx, t)=-2\nu \dfrac{D^2\varphi(\bx,t)}{\varphi(\bx,t)}+2\nu \dfrac{\nabla\varphi(x,t)\otimes \nabla\varphi(x,t)}{\varphi(\bx,t)^2}, $$
so that
$$\Delta u=\text{trace}\big (D^2u(\bx, t)\big )=-2\nu \dfrac{\Delta \varphi(\bx,t)}{\varphi(\bx,t)}+2\nu \dfrac{\Vert \nabla\varphi(x,t)\Vert^2}{\varphi(\bx,t)^2}.$$
Hence plugging this  into the Burgers equation, we have
\remi{\begin{equation*}
\begin{split}
\dpar{u}{t}+\frac{1}{2}\Vert \nabla u\Vert ^2-\nu \Delta u&=-2\nu\dfrac{\dpar{\varphi(\bx,t)}{t}}{\varphi(\bx,t)}+\frac{1}{2}\bigg ( 4\nu^2\dfrac{\Vert \nabla\varphi(\bx,t)\Vert ^2}{\varphi(\bx,t)^2}\bigg )\\
&\qquad \qquad-\nu\bigg ( -2\nu \dfrac{\Delta \varphi(\bx,t)}{\varphi(\bx,t)}+2\nu \dfrac{\Vert \nabla\varphi(x,t)\Vert^2}{\varphi(\bx,t)^2}\bigg )\\
&=-\frac{2\nu}{\varphi(\bx,t)} \bigg ( \dpar{\varphi(\bx,t)}{t}-\nu \Delta \varphi(\bx,t)\bigg )
\end{split}
\end{equation*}}
so that in the end we see that $\varphi$ needs to satisfy
$$\dpar{\varphi(\bx,t)}{t}-\nu \Delta \varphi(\bx,t)=0$$
with $\varphi=1$ on the Dirichlet  boundary and $\varphi\geq 0$ on $\Omega$.

Unfortunately, the time dependant problem
$$\dpar{u}{t}+ \Vert \nabla u\Vert^2-1=\nu \Delta u$$
does not go through as well, but this is not an issue because this is not the problem we want to solve. We want to solve
$$ \Vert \nabla u\Vert^2-1=\nu \Delta u$$ for which we again use the change of variable
$$u(\bx)=\alpha \log\big ( \varphi(\bx)\big )$$ with $\alpha$ to be determined.
We get
\begin{equation}\label{calcul}
\begin{split}
\Vert \nabla u\Vert^2-1-\nu \Delta u&=\alpha^2\dfrac{\Vert \nabla\varphi(\bx,t)\Vert ^2}{\varphi(\bx,t)^2}-1-\nu\bigg ( \alpha \dfrac{\Delta \varphi(\bx,t)}{\varphi(\bx,t)}-\alpha \dfrac{\Vert \nabla\varphi(x,t)\Vert^2}{\varphi(\bx,t)^2}\bigg )\\
&=\frac{-1}{\varphi(\bx,t)}\bigg ( \alpha\nu \Delta \varphi(\bx,t)+\varphi(\bx,t)\bigg )+\dfrac{\alpha^2+\nu\alpha}{\varphi(\bx,t)^2} \Vert \nabla\varphi(\bx,t)\Vert ^2
\end{split}
\end{equation}
so we take $\alpha=-\nu$ and we need to solve in $\Omega$
\begin{equation}
\label{eq:3}
\nu^2 \Delta \varphi(\bx,t)=\varphi(\bx,t).
\end{equation}
with the boundary condition
\begin{equation}
\label{eq:3.2}
\varphi(\bx,t)=1, \qquad \bx\in \Gamma_D.
\end{equation}
On $\Gamma_S$, inspired by  what is done for the inviscid problem, and using \eqref{calcul}, we see that a condition is
$$\nu^2 \Delta \varphi(\bx,t)\leq \varphi(\bx,t)$$ on $\Gamma_S$. This looks a bit like an obstacle problem, but this is not exactly the same (because the "obstacle" is at the boundary). In the next section, inspired by what is done for the Eikonal equation, we will propose a discretisation of this kind of condition. \remi{In the numerical section, we will also compare this boundary condition with more natural ones, such as a Neuman condition on the distance function, on $\Gamma_S$.}

\section{Numerical discretisation}
\subsection{Formulation}
We consider a triangulation of \remi{the polygonal domain} $\Omega$ that respects $\Gamma_D$ and $\Gamma_S$. They consists of triangles (or tetrahedrons) that are generically denoted by $K$.  The vertices of the triangulation are denoted by $\ba_i$, \remi{$i=1, \ldots , n_s$}. The number of element is $n_e$. For any vertex $\ba_i$, $\mathcal{V}(j)$ is the set of vertices connected to $\ba_i$ by an edge of the triangulation. Often, we make the identification between a vertex $\ba_i$ and its index $i$.  For the sake of simplicity we only consider the two dimensional case, the three dimensional one can be done in a similar way. The approximation space is 
$$\remi{V^h=\{\psi\in H^1(\Omega), \forall K, \psi_{|K}\in \P^1(K)\}\cap \{\psi=1 \text{ on }\Gamma_D\}}.$$

The trial space is
$$\remi{W^h=\{\psi\in H^1(\Omega), \forall K, \psi_{|K}\in \P^1(K)\}\cap \{\psi=0 \text{ on }\partial \Omega\}}.$$
 We write the problem as:
 Find $\varphi^h\in V^h$ such that for any $\psi^h\in W^h$,
 \begin{subequations}\label{viscous:HJ}
 \label{viscous:HJ:num}
 \begin{equation}
 \label{viscous:HJ:num:1}
 \nu^2 \int_\Omega \nabla \varphi^h\cdot \nabla \psi^h \; d\bx + \int_\Omega \varphi^h \psi^h \; d\bx =0
 \end{equation}
 coupled to boundary conditions on $\Gamma_S$. If $\Gamma_S=\emptyset$, there is nothing more to do.
 
 In the case when $\Gamma_S\neq\emptyset$, we define the boundary conditions according to what is done for the Eikonal equations, see for example \cite{abgrallHJ}. There is defined a numerical Hamiltonian $\mathcal{H}$ which role is to translate the viscosity inequality
 $$\Vert \nabla u\Vert -1-\nu^2 \Delta u\geq 0$$ on $\Gamma_S$  in the limit $\nu\rightarrow 0$. There are several possible versions, but the best (because the gradient of the numerical solution is controlled, see \cite{abgrallHJ}) is to use a Godunov Hamiltonian which amounts to write for any vertex $\ba_i\in \Gamma_S$, that
 $$\max\limits_{j\in \mathcal{V}(i)}\big ( \frac{u_i-u_j}{\Vert a_ia_j\Vert}-1\big ) =0$$
 that is the distance function $d$ satisfies
 $$u_i=\min\limits_{j\in \mathcal{V}(i)} \big ( u_j+\Vert a_ia_j\Vert\big ).$$Keeping in mind that  the solution of \eqref{viscous:HJ} is related to the solution of \eqref{viscous:HJ2} by
 $d=-\nu \log \varphi$, we will consider the following \remi{implementation } for the Soner boundary condition: for $\ba_i\in \Gamma_S$,
 \begin{equation}
 \label{viscous:HJ:num:BC}
 u_i=\exp\big (-\frac{\varphi_i}{\nu}\big ), \quad \varphi_i=\max\limits_{j\in \mathcal{V}(i)}\big ( -\nu\log u_j+\Vert a_ia_j\Vert \big ).
 \end{equation}
 \end{subequations}

\subsection{Numerical procedure}
We use the following notations. A triangle $K$ has 3 vertices, denoted by $a_i$, $a_j$, $a_k$. We assume that the elements are oriented positively. The gradient of the basis function, $\theta_i$ associated to the vertex $a_i$ is
$${\nabla \theta_i}_{|K}=\dfrac{\ba_j\ba_k^\bot}{2|K|}$$ where, for any vector $\bx$, $\bx^\bot$ is orthogonal to $\bx$ such that the basis $(\bx, \bx^\bot)$ is direct.  As usual, the angle at $\ba_i$ in $K$ is denoted by $\alpha_i^K$.

 The variational formulation in $\Omega$ leads to
\begin{subequations}
\label{discrete}
\begin{equation}\label{discrete:eq1}
M\varphi+\nu^2 R\varphi=0
\end{equation}
with the boundary condition 
\begin{equation}\label{discrete:eq2}
\varphi_i=1, \text{ for any } \ba_i\in \Gamma_D
\end{equation}
\end{subequations}
and \eqref{viscous:HJ:num:BC} on $\Gamma_S$ when this set is not empty. Here, $M$ is the mass matrix and $R$ the rigidity matrix,
$$\remi{M_{ij}=\int_\Omega \theta_i\theta_j\;d\bx, R_{ij}=\int\nabla\theta_i\cdot \nabla\theta_j\; d\bx.}$$

If $\Gamma_S=\emptyset$, this can be solved by an iterative or a direct solver. Here we have chosen the direct solver PastiX \cite {pastix}. If $\Gamma_S\neq \emptyset$ the problem becomes non linear. In that case we use an \remi{Uzawa-type} procedure: we construct a sequence of functions by initialising with $\varphi^0=1$ on $\Omega$, and from $\varphi^n$ we construct $\varphi^{n+1}$ by setting:
\begin{enumerate}
\item We compute $\widetilde{\varphi}^{n+1}$ solution of  \eqref{discrete:eq1} with the Dirichlet boundary 
\begin{equation}
\label{BC:auxiliary}\widetilde{\varphi^{n+1}}=1 \text{ on }\Gamma_D \text{ and } \widetilde{\varphi^{n+1}}=\varphi^n \text{ on } \Gamma_S.
\end{equation}
\item Then we set
\begin{equation}\label{BC:auxiliary2}
\begin{split}
\varphi^{n+1}&=\widetilde{\varphi^{n+1}} \text{ on }\Omega\backslash\Gamma_S\\
\varphi_i^{n+1}&=\exp\big (-\frac{v_i}{\nu}\big ), \quad v_i=\max\limits_{j\in \mathcal{V}(ji}\big ( -\nu\log \widetilde{\varphi^{n+1}}_j+\Vert a_ia_j\Vert \big ) \text{ on }\Gamma_S.
\end{split}
\end{equation}
\end{enumerate}

\begin{figure}[h]
\begin{center}
\includegraphics[width=0.4\textwidth]{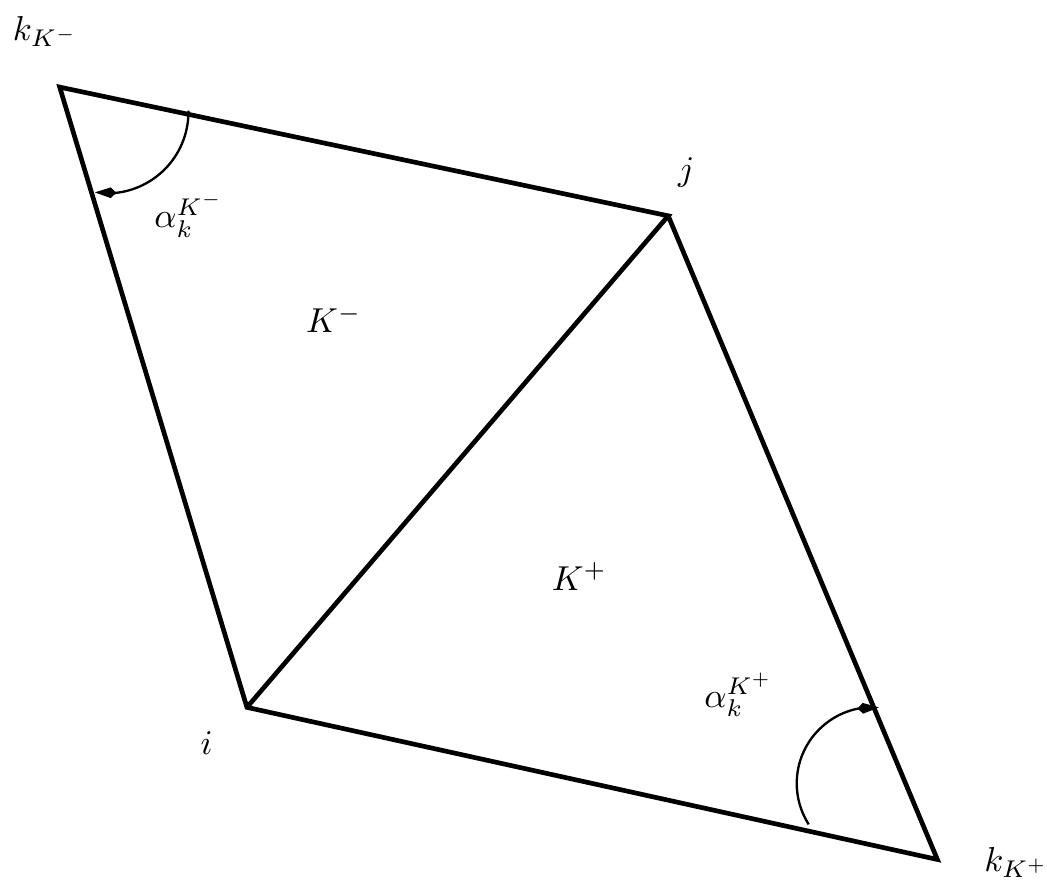}
\end{center}
\caption{\label{triangulation} Notations.}
\end{figure}

 It is well known that for $\P^1$ approximation on triangular elements, 
the contribution of $K$ for $R$ is
$$\frac{1}{2}\begin{pmatrix}
\cot\alpha_j^K+\cot\alpha_k^K & -\cot\alpha_k^K & -\cot\alpha_j^K\\
-\cot\alpha_k^K & \cot\alpha_i^K+\cot\alpha_k^K & -\cot\alpha_i^K\\
-\cot\alpha_j^K & -\cot\alpha_i^K & \cot\alpha_i^K+\cot\alpha_j^K
\end{pmatrix}
$$
Since $\cot\alpha+\cot\beta=\frac{\sin(\alpha+\beta)}{\sin\alpha\sin \beta}$, and since $\alpha_i^K+\alpha_j^K+\alpha_k^K=\pi$, the diagonal terms are positive \footnote{a quicker way to see this is to write $R_{ii}=\frac{\Vert \ba_j-\ba_k\Vert^2}{2|K|}$.}. The term $R_{ij}$ for two adjacent points is, since $\ba_i$ and $\ba_j$ defines the common edges between two triangles $K^+$ and $K^-$,
$$\remi{R_{ij}=\int_{K^+}\langle \nabla\theta_i, \nabla\theta_j\rangle +\int_{K^-}\langle \nabla\theta_i, \nabla\theta_j\rangle=-\frac{1}{2}\big (\cot\alpha_k^{K^+}+\cot\alpha_k^{K^-}\big )}$$
and it is known that  $R_{ij}\leq 0$ for $i\neq j$ if and only if 
\begin{equation}
\label{positivity}\alpha_k^{K^+}+\alpha_k^{K^-}\leq \pi,
\end{equation}see figure \ref{triangulation}.

The mass matrix $M$ is a matrix with positive entries, and the contribution $M_K$  of element $K$ to it is
$$|K|\begin{pmatrix} 1/6 & 1/12& 1/12\\ 1/12& 1/6 &1/6\\  1/12 & 1/12 & 1/6 \end{pmatrix}$$

If for any vertex, the condition \eqref{positivity} is met, the solution of the auxiliary problem satisfies $$\remi{\min\limits_{\bx\in \Gamma_S}\varphi^n(\bx)\leq \widetilde{\varphi}^{n+1}\leq 1},$$
and from \eqref{BC:auxiliary2}, we see that $$0\leq \min\limits_{\bx\in \Gamma_S}\varphi^{n+1}(\bx)\leq {\varphi}^{n+1}\leq 1.$$
Similarly, we can show that the sequence is monotone increasing, and since $\varphi^n\geq 0$, it is convergent. The monotone increasing nature comes from $\varphi^1\leq 1=\varphi^0$ and then we proceed by induction. We have (using the discrete maximum principle thanks to the condition \eqref{positivity}) that $\widetilde{\varphi^{n+1}}\geq \varphi^n$ and then using \eqref{BC:auxiliary2}, we see that $\varphi^{n+1}\leq \widetilde{\varphi^{n+1}}$ on $\Gamma_S$.
We have thus shown:
\begin{proposition}
If the variational formulation of the Laplace operator satisfies a maximum principle, the sequence $(\varphi^n)_{n\in \N}$ converges. This is in particular true is the triangulation satisfies the angle condition
\eqref{positivity}.
\end{proposition}

\section{Numerical examples}
All the calculations have been done on an Imac with 3.5GHz Quad-Core Intel core i7 processors with the version 6.0.2 of PastiX \cite{pastix}. To report the performance of the solver, for a mesh with 295  296vertices, 587 520 elements generated by GMSH \cite{gmsh} using the frontal Delaunay option, the symbolic factorisation takes 0.47 s, the evaluation of the non zeros entries of the matrices and the \remi{second} hand side takes 0.11 s, and the solution takes 9.16 seconds. The averaged maximal band-with of the matrix was 171 781, its maximal band-with is 292 530. The computations have been done sequential. We do only calculations when $\Gamma_S\neq \emptyset$ because they are a priori more complicated. \remi{The symbolic factorisation is done once for all (for a given mesh) if the Uzawa-type method is needed.}

\begin{figure}[h]
\begin{center}
\subfigure[]{\includegraphics[width=0.45\textwidth]{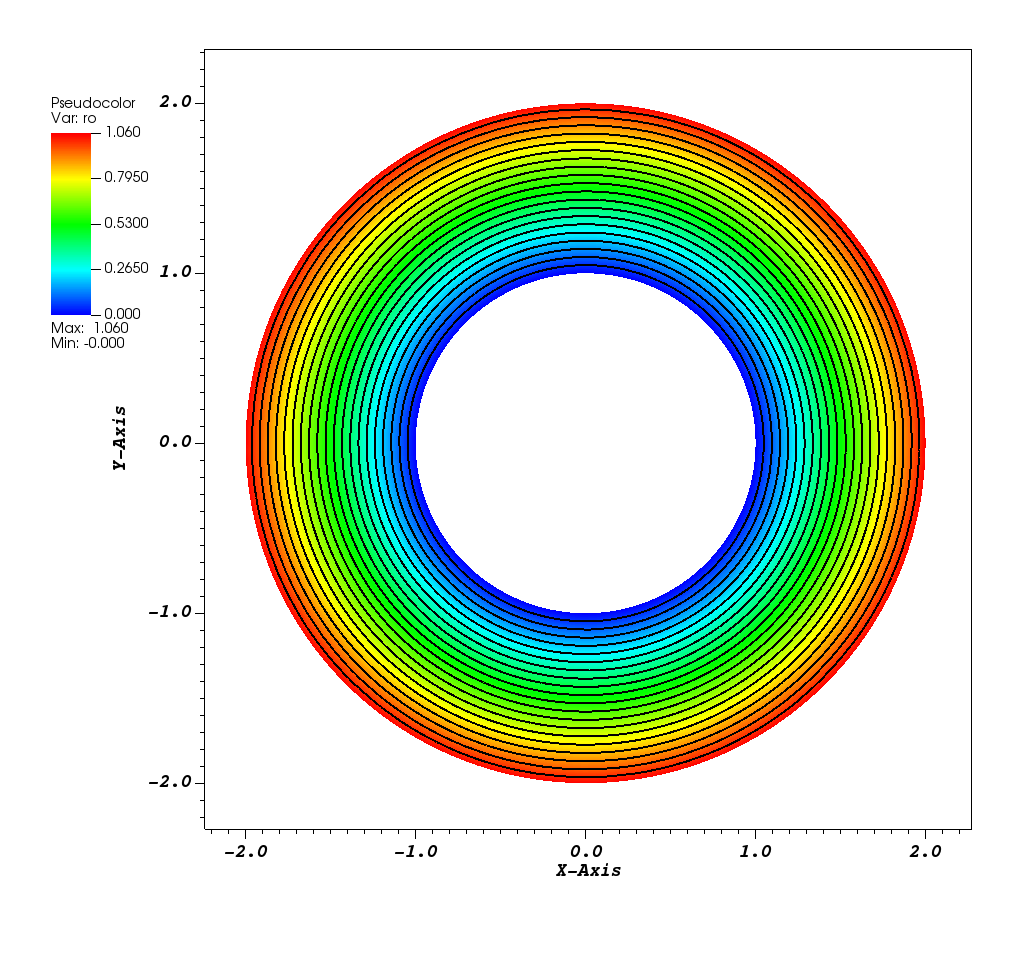}}
\subfigure[]{\includegraphics[width=0.45\textwidth]{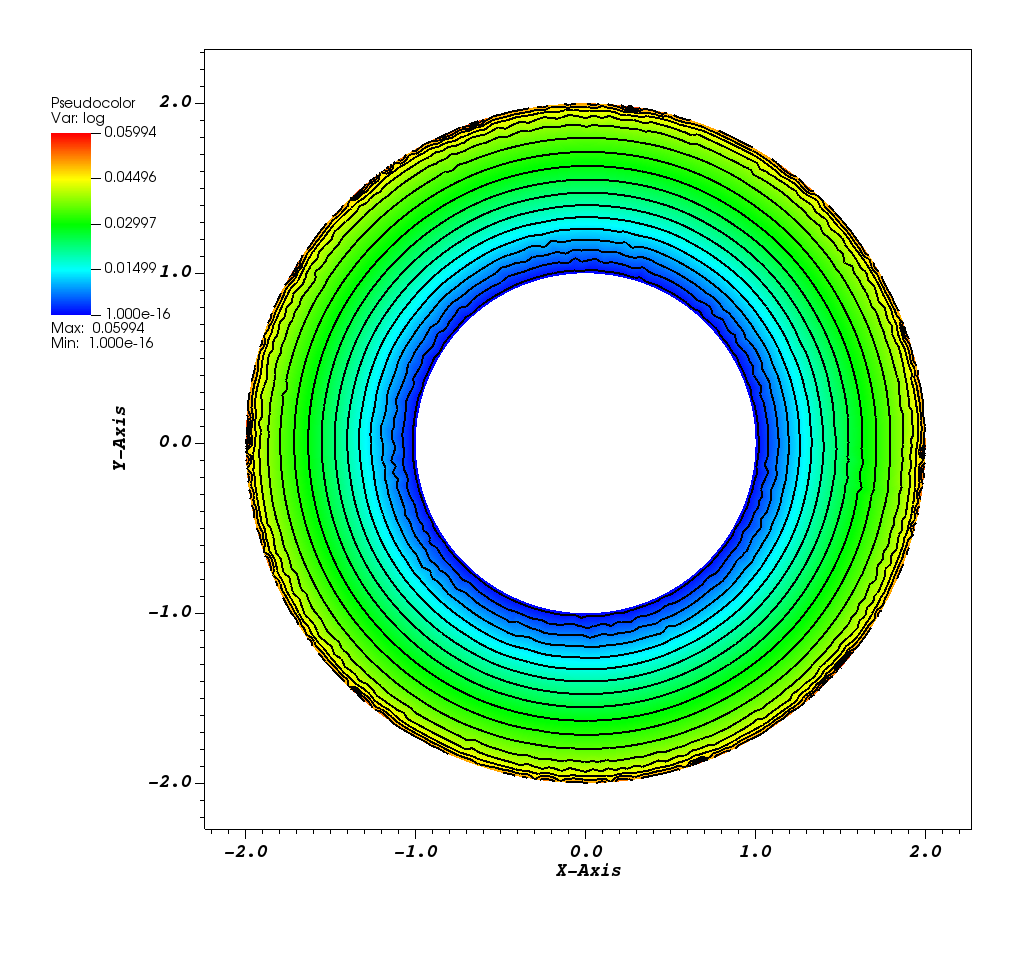}}
\end{center}
\caption{\label{annulus} Results for the distance in an annulus. On (a) we have the solution, and on (b) we have $\vert -\nu \log u_i-d_i\vert$}
\end{figure}

The first test is the evaluation of the distance function  on the annulus $\{\bx, 1\leq \Vert\bx\Vert\leq 2\}$ The viscosity is set to $\nu=0.1$. \remi{The Dirichlet condition is set for the inner circle, and the Soner condition on the outer circle.} The solution is displayed on Figure \ref{annulus}-a, while the error to the true solution is on Figure \ref{annulus}-b.

\begin{figure}[h]
\begin{center}
\subfigure[]{\includegraphics[width=0.45\textwidth]{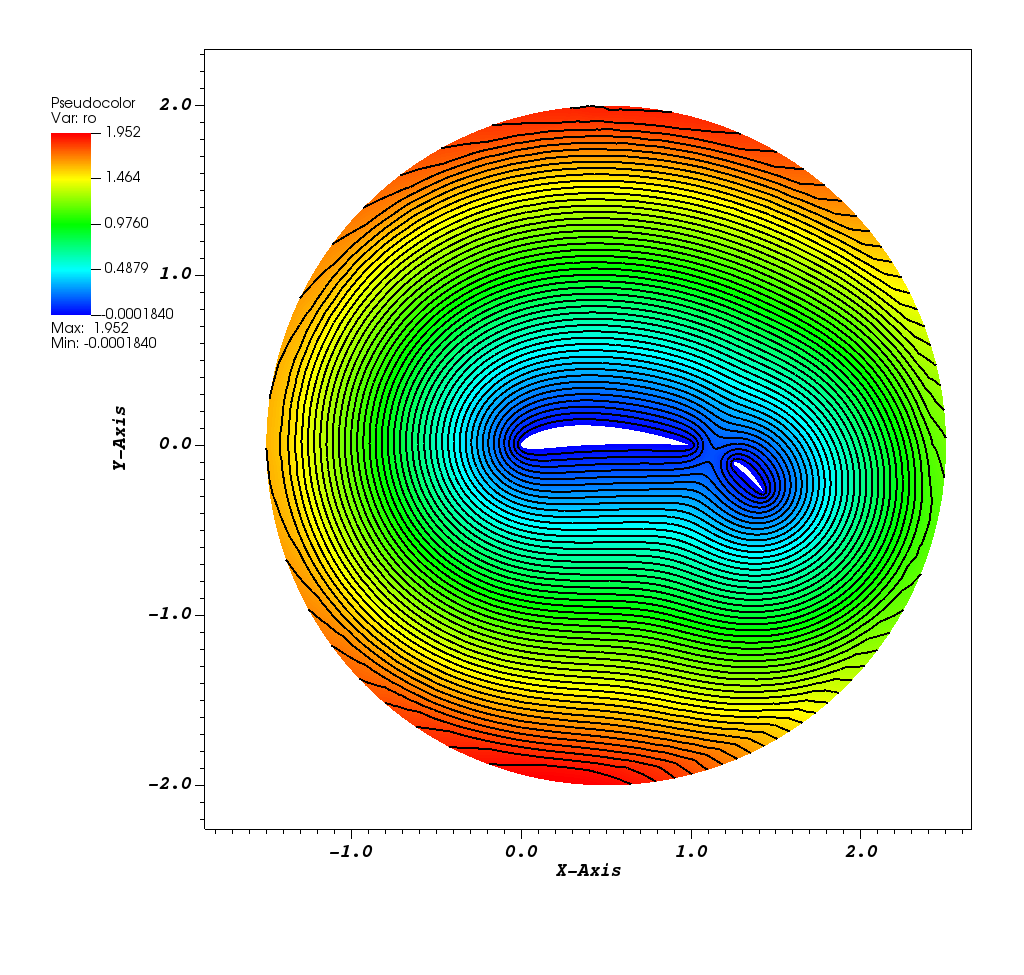}}
\subfigure[]{\includegraphics[width=0.45\textwidth]{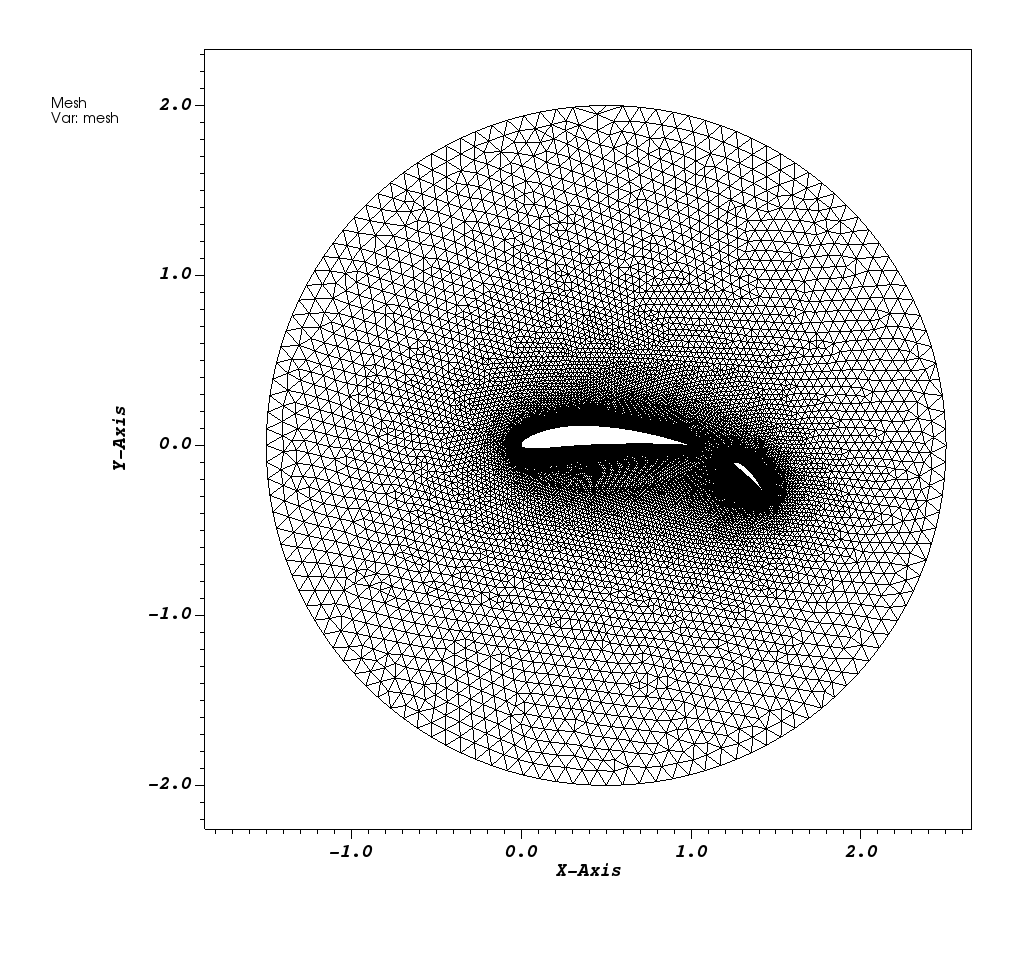}}
\subfigure[]{\includegraphics[width=0.45\textwidth]{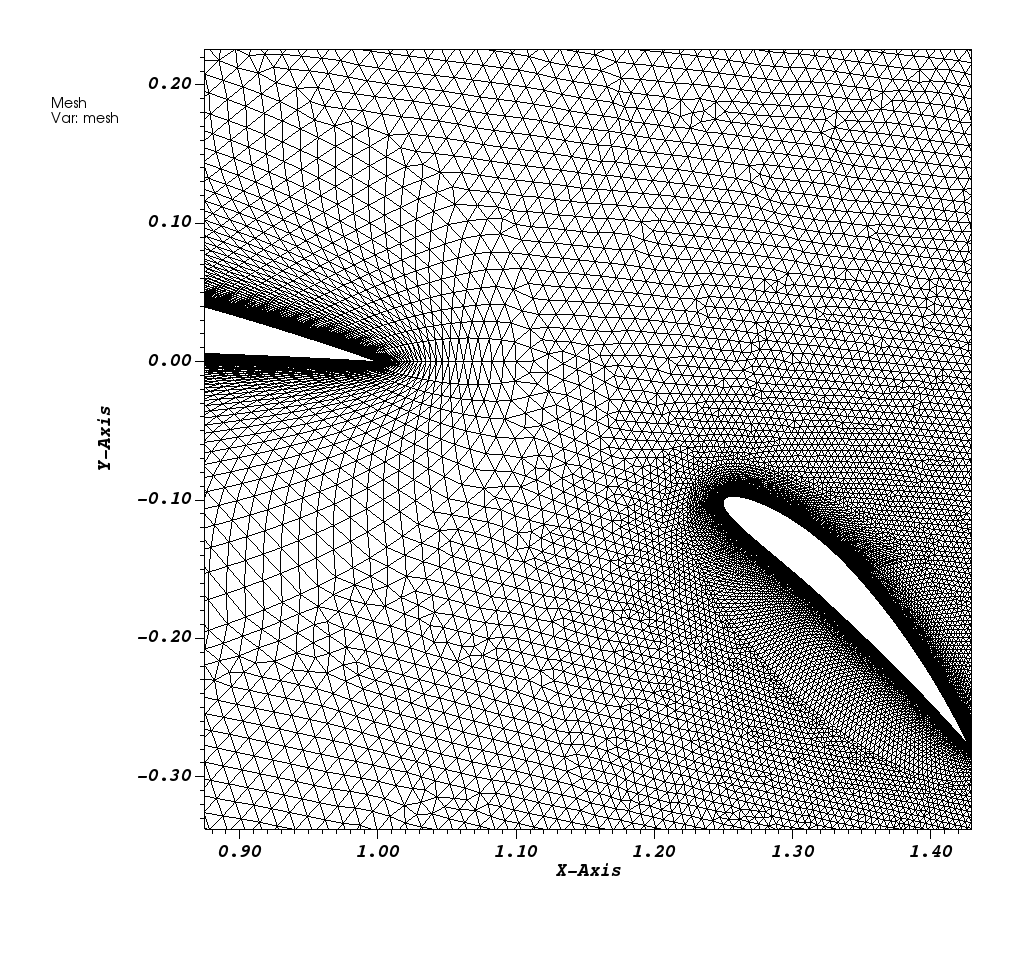}}
\end{center}
\caption{\label{Naca} Results for the two Naca problem. The mesh and a zoom of the mesh is displayed.}
\end{figure}

The second case is that of the distance to a body made of two NACA airfoils. \remi{The Dirichlet condition is set on the airfoils, and the Soner one on the outer boundary.} The results are displayed on figure \ref{Naca} as well as the mesh. Here again, $\nu=0.1$

\remi{We also have considered a case where $\Gamma_D$ and $\Gamma_S$ are not disjoint. The example under consideration is
\begin{subequations}\label{eq:aslam}
\begin{equation}\label{eq:aslam:1}\Omega=\{\bx=(x_1,x_2) \in \R^2, 1\leq \Vert \bx\Vert \leq 2\Remi{,}  x_1\geq 0\}.\end{equation}
We take \begin{equation}\label{eq:aslam:2}\Gamma_D=\{(x_1,0), 1\leq x_1\leq 2\}\text{ , and }\Gamma_S=\partial \Omega\backslash \Gamma_S.
\end{equation}
\end{subequations}
We take $\nu=0.01$.
We provide the result on a fine mesh ($151\;713$ points and       $301\;568$ elements) on figure \ref{aslam}.
\begin{figure}[h]
\begin{center}
\includegraphics[width=0.8\textwidth]{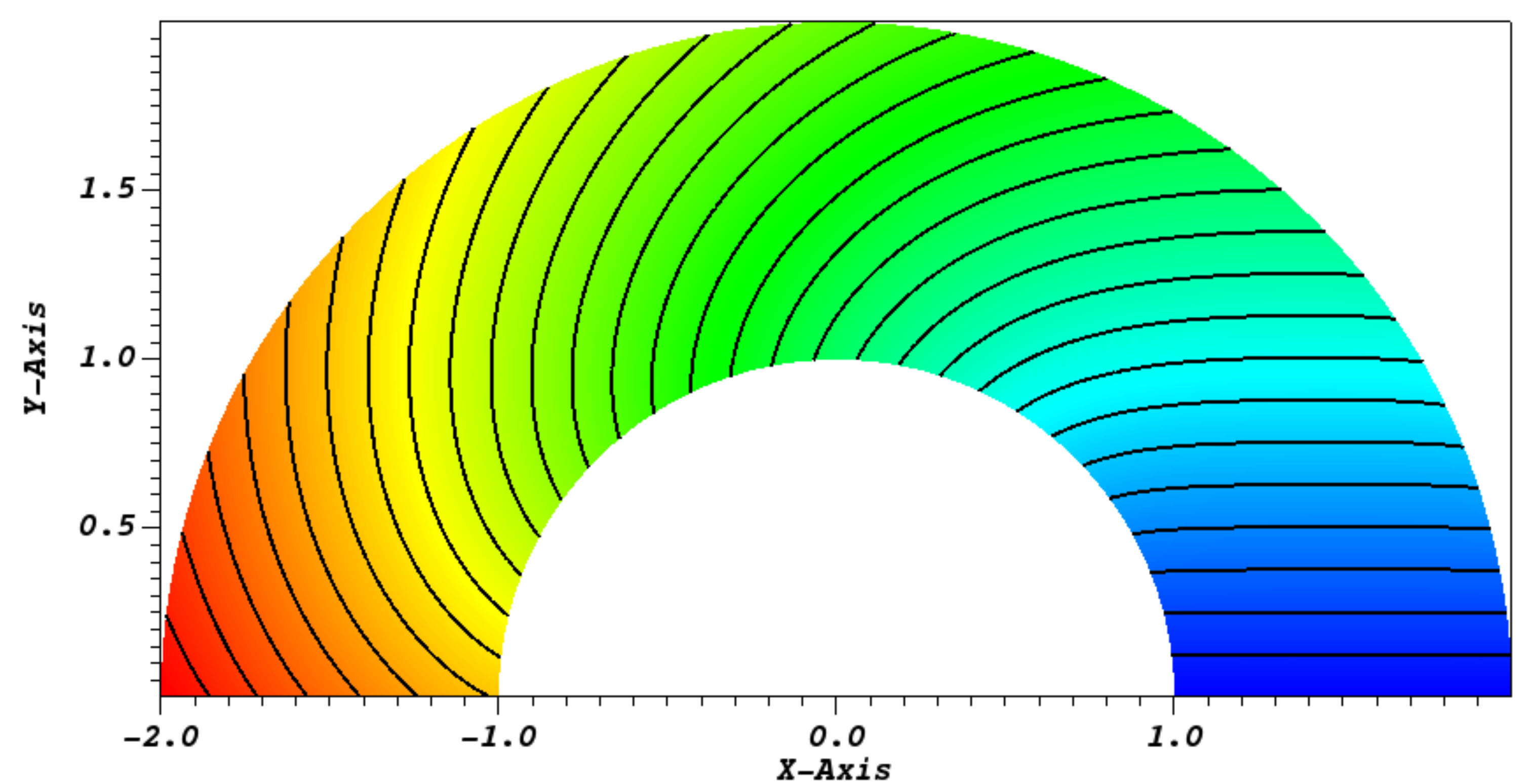}
\end{center}
\caption{\label{aslam} Solution for the geometry and the boundary conditions defined by \eqref{eq:aslam}. The values range between $0$ and $3.907$, 30 isolines.}
\end{figure}

In figure \ref{convergenceiterative}, we compare the iterative convergence of the algorithm for several meshes ($2473$, $9657$ and  $151\;713$ vertices).
\begin{figure}[h]
\begin{center}
\includegraphics[width=0.8\textwidth]{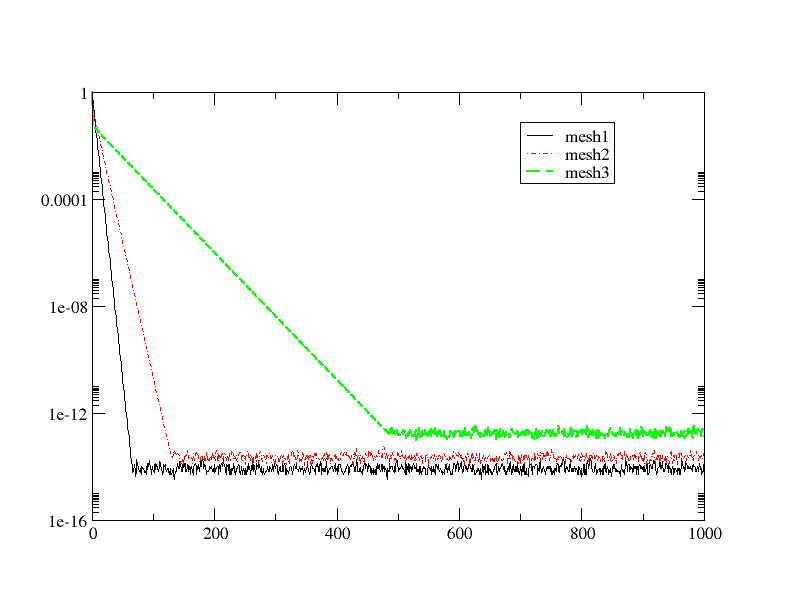}
\caption{\label{convergenceiterative} Convergence to the steady solution for several meshes.}
\end{center}
\end{figure}
We observe that we need 65 iterations for the coarse mesh, 130 for the medium one  and 485 for the fine one. The \Remi{ratio} of mesh points w.r.t. the coarse one is 1:4:61, and for the iteration the ratios  1:2:7.5, so we see that the cost evolves like $h^{-1}$ (the mesh is very regular). For a classical explicit hyperbolic solver, the cost scales like the number of points, so $h^{-2}$. In this case, and others that we have computed (such as the airfoil case, the conclusion is similar, or even better. }

\remi{To end this paragraph, let us comment a bit on the boundary condition of Soner type.
Since we expect that approximation of the distance has a gradient of norm approximately equal to unity, one may wonder why a Neuman type boundary condition on the distance, say $\nabla u\cdot \bn=1$, would not fit. Indeed, this has been our starting point for imposing boundary conditions on $\Gamma_S$. 

If we do that, we first have, since $u=-\nu \log \varphi$, 
$$\dpar{u}{n}=\nabla u\cdot \bn=-\nu \dfrac{\nabla \varphi \cdot\bn}{\varphi}.$$
Thus, to set $\dpar{u}{n}=1$,  we are let to Robin type conditions on $\varphi$. We have compared our formulation and this one on the fine mesh of the previous test case to the exact solution
$$
u(x,y)=\left \{\begin{array}{ll}
y & \text{ if } 1\leq x\leq 2 \text{ and } (x,y)\in \Omega\\
d_1+\theta & \text{ else,}
\end{array}\right .
$$
where, if $\bx=(x,y)$,
$$d_1=\Vert \bx-\mathbf{p}\Vert \text{ with }\mathbf{p}=\frac{\bx}{\Vert\bx\Vert^2}-\sqrt{1-\frac{1}{\Vert \bx\Vert^2}} \frac{\bx^\bot}{\Vert \bx\Vert}$$ and $\theta$ is the arclength on the inner circle between $\mathbf{p}$ and $(1,0)$.
The results are displayed on figure \ref{comparison} with $\nu=0.01$. 
\begin{figure}[h]
\begin{center}
\subfigure[]{\includegraphics[width=0.45\textwidth]{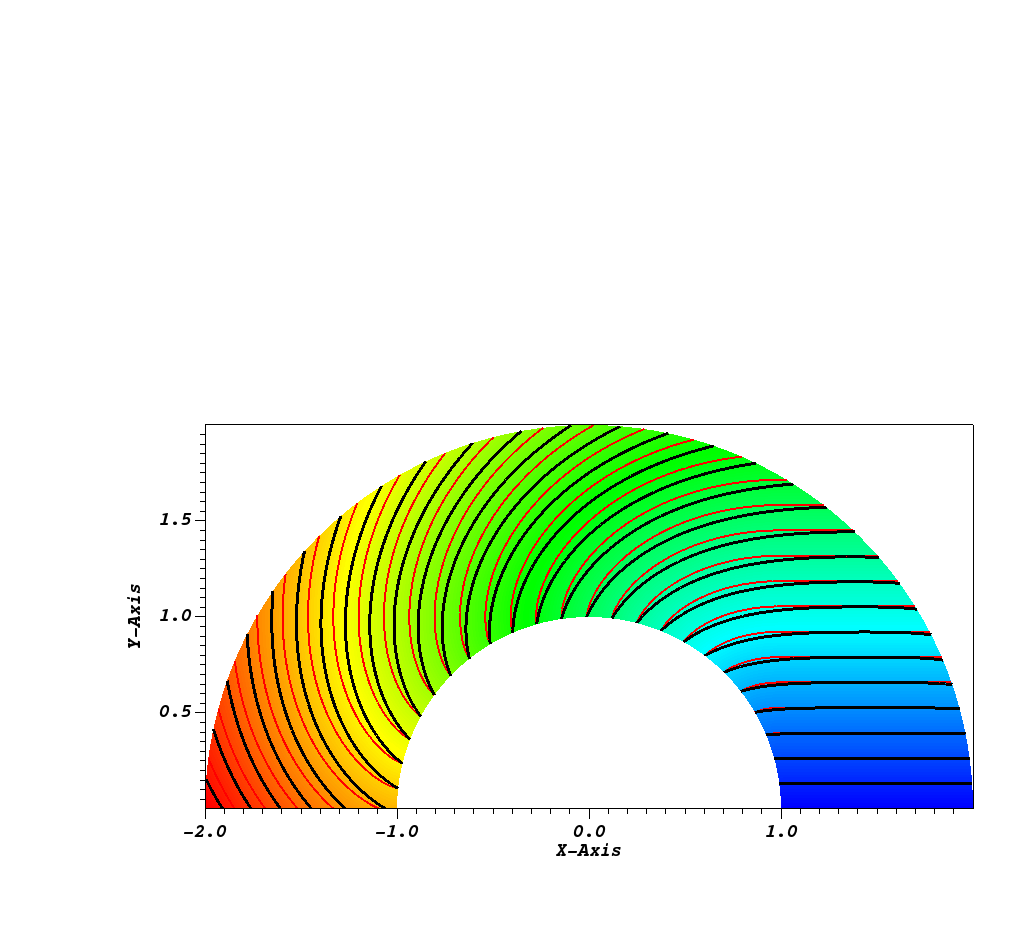}}
\subfigure[]{\includegraphics[width=0.45\textwidth]{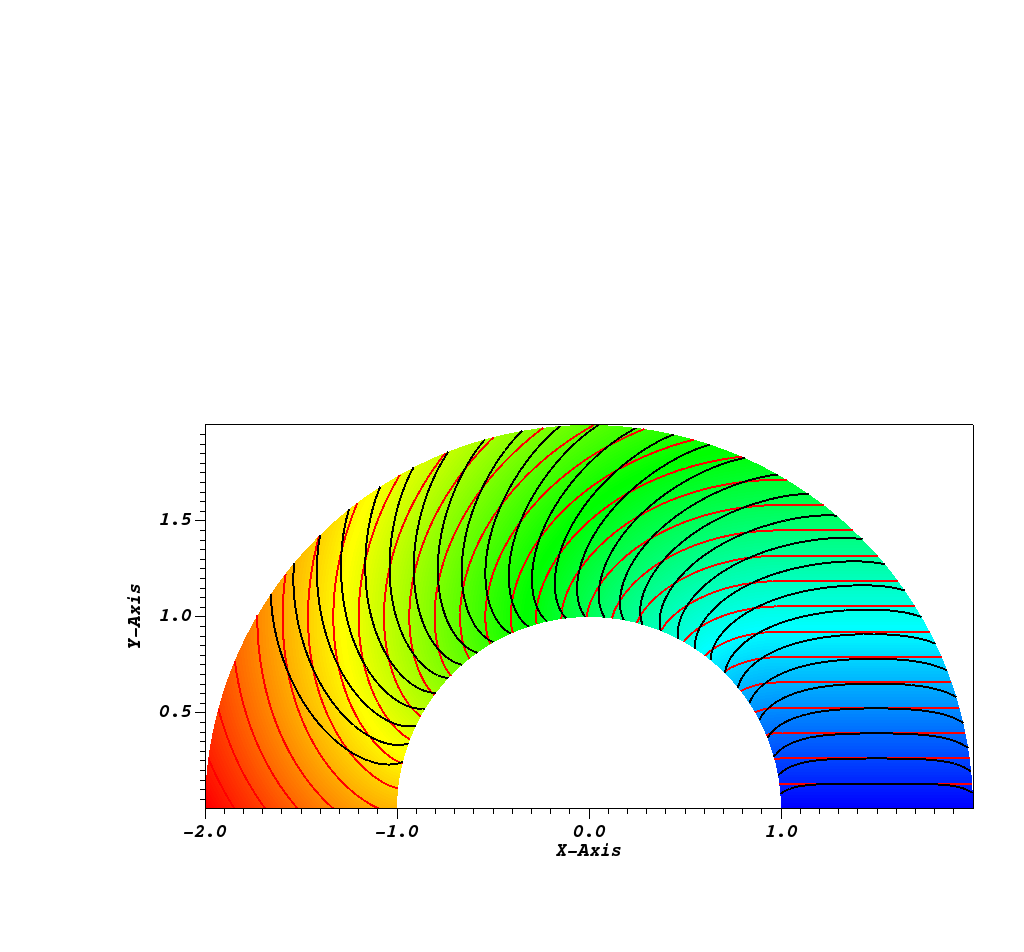}}
\end{center}
\caption{\label{comparison} (a): isolines for the exact solution (red) and the solution with the Hopf-Cole transform, (b): isolines for the exact solution(red) and the solution with Robin boundary conditions (black). The background colour is that of the exact solution, in both case we have 30 isolines between $0$ and $u=3.826$ (the maximum value of the distance on this domain). The solution with the Neumann condition get values larger than 3.9, this explains why the isolines stop}
\end{figure}
From the figure it is clear that the solution with the Neumann condition is completely off compared with our method, though this one is not perfect (because of $\nu$ not small enough).
However, the solution with the Robin condition could be used as an initial guess, this has not been done in our code.}

\section{Conclusion}
This paper is dedicated to Roland Glowinski. He always have been very nice to the author. Roland has also worked a lot on Hamilton Jacobi equations, for example in \cite{zbMATH07161449,zbMATH07074516,zbMATH06993413}. This plus the handling  of \cite{nishi}  as editor of JCP has motivated the present paper. 

A final remark is that it is certainly possible to establish rigorous error bounds between the distance function and what is computed here using approximation results between the viscous regularisation of Hamilton Jacobi equations and standard $L^\infty$ error estimates on $\P^1$ approximation of elliptic equations. This \remi{has not been} done here because we were motivated by designing a working algorithm.

This algorithm has its own drawbacks. The first one is that when $\nu$ becomes very small, the problem becomes stiffer and stiffer. When the domain is large, the actual value of the solution of the elliptic problem becomes extremely small. It is interesting to note links with large deviation problems (see the last chapter of \cite{barles} where exactly the same PDE is studied, for completely different reasons).
However, if one comes back to the initial motivation of this work (finding the distance function for turbulence modelling), our experience is that the computation close to  the Dirichlet boundary is very reliable, and that applying the Dirichlet condition on all boundaries is enough to get a good approximation. 

When the mesh is \remi{too distorted} so that a discrete maximum principle does not apply, the solution  can be slightly above $1$ (so that the 'distance' would be negative): \remi{in that case, }the solution provided by this method can be used as a good initial condition to a 'traditional ' Hamilton Jacobi problem.

\bibliographystyle{plain}
\bibliography{biblio}
\end{document}

%% file: command.tex
\usepackage{fullpage}

\usepackage{amsmath,amsfonts,amsthm,amssymb}
\usepackage{algorithm,algorithmic}
\usepackage{mathrsfs}
\usepackage{color,xcolor}
\usepackage{ifpdf}
\usepackage{psfrag}
\usepackage{graphicx,graphics,subfigure,epsfig}
\usepackage{url}
\usepackage{hyperref}
\usepackage{xspace}
\usepackage{xargs}
\usepackage[colorinlistoftodos,prependcaption,textsize=tiny]{todonotes}
 \usepackage{ booktabs,array}

\newtheorem{theorem}{Theorem}[section]

\newtheorem{proposition}[theorem]{Proposition}

\definecolor{OliveGreen}{rgb}{0.33, 0.42, 0.18}
\definecolor{Plum}{rgb}{0.8, 0.6, 0.8}
\newcommandx{\unsure}[2][1=]{\todo[linecolor=red,backgroundcolor=red!25,bordercolor=red,#1]{#2}}
\newcommandx{\change}[2][1=]{\todo[linecolor=blue,backgroundcolor=blue!25,bordercolor=blue,#1]{#2}}
\newcommandx{\info}[2][1=]{\todo[linecolor=OliveGreen,backgroundcolor=OliveGreen!25,bordercolor=OliveGreen,#1]{#2}}
\newcommandx{\improvement}[2][1=]{\todo[linecolor=Plum,backgroundcolor=Plum!25,bordercolor=Plum,#1]{#2}}
\newcommandx{\thiswillnotshow}[2][1=]{\todo[disable,#1]{#2}}
\newcommand{\dpar}[2]{\dfrac{\partial #1}{\partial #2}}

 \newcommand{\R}{\mathbb R}
 
 \newcommand{\N}{\mathbb N}

\renewcommand{\P}{\mathbb P}

\newcommand{\bn}{{\mathbf {n}}}

\newcommand{\ba}{\mathbf{a}}
\newcommand{\bx}{\mathbf{x}}

\newcommand{\remi}[1]{#1}

\hypersetup{
pdftitle={A general framework to construct schemes that satisfy additional conservation relations.\\
Application to entropy conservative and entropy dissipative schemes}
pdfauthor={R. Abgrall},
pdfpagemode=UseOutlines,
linkbordercolor=0 0 0,
linkcolor=red,
citecolor=blue,
colorlinks=true,
bookmarks = true
}